\documentclass[12pt]{amsart}
\usepackage[]{amsmath, amsthm, amsfonts, amssymb, epsfig}
\newcommand\C{{\mathbb C}}

\newcommand\dee{\partial}
\newcommand\Om{\Omega}
\newcommand\Obar{\overline{\Omega}}
\newcommand\Ot{\widetilde\Omega}
\newcommand\Ohat{\widehat\Omega}

\renewcommand\phi{\varphi}

\numberwithin{equation}{section}

\begin{document}

\title[Hejhal's theorem]{Ruminations on Hejhal's theorem about the
Bergman and Szeg\H o kernels}
\author[S.~Bell \& B.~Gustafsson]
{Steven R.~Bell* and Bj\"orn Gustafsson}

\address[]{Mathematics Department, Purdue University, West Lafayette,
IN  47907}
\email{bell@purdue.edu}

\address[]{Department of Mathematics, KTH, 100 44 Stockholm, Sweden}
\email{gbjorn@kth.se}



\begin{abstract}
We give a new proof of Dennis Hejhal's theorem on
the nondegeneracy of the matrix that appears in the
identity relating the Bergman and Szeg\H o kernels
of a smoothly bounded finitely connected domain in
the plane. Mergelyan's theorem is at the heart of
the argument. We explore connections of Hejhal's
theorem to properties of the zeroes of the Szeg\H o
kernel and propose some ideas to better understand
Hejhal's original theorem.
\end{abstract}

\maketitle

\theoremstyle{plain}

\newtheorem {thm}{Theorem}[section]
\newtheorem {lem}[thm]{Lemma}

\hyphenation{bi-hol-o-mor-phic}
\hyphenation{hol-o-mor-phic}

\begin{quotation}
\begin{center}
\begin{em}
To celebrate the legacy of Harold S.~Shapiro
\end{em}
\end{center}
\end{quotation}

\section{Introduction}
\label{sec1}

Dennis Hejhal, in a remarkable {\it tour de force\/}
that filled a volume of the Memoirs of the AMS \cite{Hej},
proved, among many other things, that the matrix of
coefficients $[\lambda_{ij}]$ that appears in the identity
relating the Bergman kernel $K(z,w)$ to the Szeg\H o kernel
$S(z,w)$,
\begin{equation}
\label{BS}
K(z,w)=4\pi S(z,w)^2 + \sum_{i,j=1}^{n-1}
\lambda_{ij}F_i'(z)\overline{F_j'(w)},
\end{equation}
in a bounded $n$-connected smoothly bounded domain in
the plane is nondegenerate and, in fact, positive definite.
Hejhal's proof of this result used a great deal of machinery
from analysis and geometry, including key use of theta functions
on Riemann surfaces. The purpose of this paper is to give
a rather short proof of the nondegeneracy of the matrix that
uses only Mergelyan's theorem and basic properties
of the Bergman and Szeg\H o kernel functions. We also
explore how these results are connected to properties
of the zeroes of the Szeg\H o kernel.

The authors stumbled upon this application of Mergelyan's
theorem after their work on double quadrature domains \cite{BGS},
which turns out to be a subject closely connected to Hejhal's
theorem. This work sprouted from the influential works of
Harold Shapiro and his many collaborators, including
\cite{AS}, \cite{S}, and \cite{SU}. Both of us have
greatly benefited from Harold's mentorship and generosity 
and so it seems fitting to offer this work in a volume in
his honor. (We must also mention here that Avci's Stanford
Thesis \cite{Av} also played an important role in our studies.)

We tried to give a new proof also of Hejhal's
full result, that the matrix is {\it positive definite}, but could achieve this only 
in a few special cases (connectivity two and three).
In our attempts we however learned a great deal
about alternate arguments and how Hejhal's
result fits in the grand scheme of things. We could not
resist including some of these observations here, giving
this paper an expository component. It has been a
great subject to ruminate upon!

\section{Preliminaries}
\label{sec2}

The transformation identities for the Bergman and
Szeg\H o kernels and the harmonic measure functions
yield that the nondegeneracy of the
matrix $[\lambda_{ij}]$ in equation \eqref{BS}
is invariant under conformal changes of variables.
Hence, we may always suppose that the domain $\Om$
under study is a bounded domain in $\C$ bounded by
$n>1$ nonintersecting smooth real analytic curves.

We now set up some definitions and notation that
we will use throughout the paper.

We denote the boundary of $\Om$ by $b\Om$ and
provide it with the standard orientation.
Let $\gamma_n$ denote the outer boundary curve of
$\Om$, and denote the inner boundary curves by
$\gamma_j$, $j=1,\dots,n-1$.

The double $\Ohat$ of $\Om$ is a compact Riemann
surface of genus $n-1$ obtained by using the Schwarz
reflection principle to glue a copy $\Ot$ of $\Om$
to $\Om$ along the boundary of $\Om$, using the
function $z$ as a chart on $\Om$ and $\bar z$ as a
chart on $\Ot$.

We now define curves that go around the $n-1$
handles of $\Ohat$. Let $\sigma_j$ be a curve in
$\Om$ that starts on the outer boundary $\gamma_n$
and ends on $\gamma_j$ for $j=1,\dots, n-1$. The
curves $\sigma_j$ can be defined so that their
closures do not intersect. Note that, in this case,
$\Om-\cup_{j=1}^{n-1}\sigma_j$ is a simply connected
domain. Let $\beta_j$ denote the curve on $\Ohat$
that first follows $\sigma_j$ in $\Om$, and then
follows the copy of $-\sigma_j$ in $\Ot$ to
connect back to the starting point. We think of
$\beta_j$ as going around the $j$-th handle of
$\Ohat$ and we note that the $n-1$ curves
$\gamma_j$, $j=1,\dots,n-1$, together with the
$n-1$ curves $\beta_j$ form a homology basis for
the double.

The Bergman kernel $K(z,w)$ is the kernel for
the orthogonal projection of $L^2(\Om)$ onto
its closed subspace of holomorphic functions
in $L^2$. The Szeg\H o kernel $S(z,w)$ is the
kernel for the orthogonal projection of $L^2(b\Om)$
onto its subspace consisting of $L^2$ boundary
values of holomorphic functions.  We refer the
reader to the classic books \cite{{Berg},{G},{N}}
for the basic facts about these kernels and to
\cite{B} for a treatment of the subject very
much in line with the approach of this paper.
In fact, this paper fills in a missing chapter
of \cite{B}.

The functions $F_j'(z)$ appearing in equation~\eqref{BS}
are given by
\begin{equation}
\label{Fw}
F_j'(z)=2\frac{\dee \omega_j}{\dee z},
\end{equation} 
where $\omega_j$
is the harmonic function on $\Om$ that has boundary
values equal to one on $\gamma_j$ and equal to zero
on the other boundary curves. The notation is
traditional; $F_j'$ is locally the derivative of
the holomorphic function with real part $\omega_j$,
but it is not globally the derivative of a holomorphic
function on $\Om$.

We let $\Lambda(z,w)$ denote the complimentary kernel
to the Bergman kernel which satisfies
\begin{equation}
\label{BL}
K(z,w)\,dz=-\overline{\Lambda(z,w)\,dz}
\end{equation}
for $z$ in $b\Om$ and $w\in\Om$. (Our choice of
symbols for the kernel functions follows \cite{B}.
In the literature, the kernels are often denoted
by only $K$ and $L$ with various tildes or
hats.)
The identity~\eqref{BL} yields that the holomorphic
one-form $K(z,w)\,dz$ on $\Om$ extends to the double as
a meromorphic one-form $\kappa_w$ by setting it equal to the
conjugate of $-\Lambda(z,w)\,dz$ on the backside of $\Om$
in the double and using the identity to connect the
definitions at the boundary.
Let $G(z,w)$ denote the classical Green's function
(with singular behavior $-\ln|z-w|$ near $z=w$).
Since $K$ and $\Lambda$
are related to the Green's function via
\begin{align}
K(z,w) &= -\frac{2}{\pi}\frac{\dee^2}{\dee z\dee\bar w}G(z,w) \label{BLG} \\
\Lambda(z,w) &= -\frac{2}{\pi}\frac{\dee^2}{\dee z\dee w}G(z,w), \nonumber
\end{align}
it follows that the periods of $\kappa_w$ about each $\beta_j$
vanish (if $w$ does not fall on any of the $\sigma_j$), i.e.,
\begin{equation}
\label{zperiods}
\int_{\beta_j} \kappa_w =0
\end{equation}
for $j=1,\dots,n-1$. This very important fact, due to Schiffer and
Spencer \cite{SSp}, will be an essential ingredient in the proof
in the next section. We briefly explain the result here to make
this paper self contained. Using the definition $\frac{\dee}{\dee z}
=\frac{1}{2}\left(\frac{\dee}{\dee x}-i \frac{\dee}{\dee y}\right)$
and writing out $\int_{\sigma_j}\frac{\dee G}{\dee z}\ dz$ yields
that the real part of the integral is given by
$\frac{1}{2}\int_{\sigma_j}\frac{\dee G}{\dee x}\,dx
+\frac{\dee G}{\dee y}\,dy$, and so the integral is pure imaginary
because $G$ vanishes at the endpoints of $\sigma_j$, which fall
on the boundary of $\Om$. Since $G$ is real valued,
the conjugate of $\dee G/\dee z$ is equal to
$\dee G/\dee\bar z$, and the conjugate of the integral is equal to
$\int_{\sigma_j}\frac{\dee G}{\dee \bar z}\ d\bar z$. Hence,
$$0=\int_{\sigma_j}\frac{\dee G}{\dee z}\ dz
+\int_{\sigma_j}\frac{\dee G}{\dee \bar z}\ d\bar z,$$
and multiplying by $-2/\pi$, differentiating with respect to $\bar w$,
and using (\ref{BL}) yields the result. (The three minus signs, one
from the conjugate of a pure imaginary integral, one from (\ref{BL}),
and one from the opposite direction of the curve, guarantee that
the integrals cancel.)

We refer the reader to the standard references for the
basic properties of $K$ and $\Lambda$. We only note here
that $K(z,w)$ is holomorphic in $z$ and antiholomorphic in
$w$, $\Lambda(z,w)$ is holomorphic in both variables off
the diagonal, $K(w,z)=\overline{K(z,w)}$,
$\Lambda(z,w)=\Lambda(w,z)$, $\Lambda(z,w)$
has a double pole in $z$ at $w$ with principal part
$\frac{1}{\pi}(z-w)^{-2}$. Both $K(z,w)$ and $\Lambda(z,w)$
extend holomorphically past the boundary in $z$ for fixed
$w$ in $\Om$, $K(z,w)$ is $C^\infty$-smooth on
$\Obar\times\Obar$ minus the diagonal $\{(z,z):z\in b\Om\}$
and $\Lambda(z,w)$ is in $C^\infty$ on $\Obar\times\Obar$
minus the diagonal $\{(z,z):z\in\Obar\}$.

The Garabedian kernel $L(z,w)$ is the complimentary kernel
to the Szeg\H o kernel and satisfies the identity
\begin{equation}
\label{SL}
\overline{S(z,w)}ds_z=\frac{1}{i}L(z,w)\,dz
\end{equation}
for $z$ in $b\Om$ and $w\in\Om$, where $ds_z$ represents the
element of arc length on the boundary.
Squaring this formula yields that
\begin{equation}
\label{SL2}
\overline{S(z,w)^2}\,d\bar z=-L(z,w)^2\,dz
\end{equation}
for $z$ in $b\Om$ and $w\in\Om$ and this shows that
the holomorphic one-form $S(z,w)^2\,dz$ extends to be
a meromorphic one-form $\sigma_w$ on the double by
defining it to be the conjugate of $-L(z,w)^2 \,dz$
on the back side of $\Om$ in the double. The key assertion for the 
paper is that we can take linear combinations of $\sigma_w$ that have
$\beta$-periods being anything we like, and this will imply the
non-degeneracy of the $\lambda$-matrix. To be precise, we have
the following theorem, to be proved in Section~\ref{sec3}.

\begin{thm} 
\label{thm:main}
The linear span of 
$$
\left\{\left(\int_{\beta_1}\sigma_{w}, \dots,
\int_{\beta_{n-1}}\sigma_{w}\right): w\in\Omega\right\}
$$
is dense in $\C^{n-1}$. As a consequence, the matrix $[\lambda_{ij}]$
is non-singular.
\end{thm}

To continue describing background material, we note that $L(z,w)=-L(w,z)$ and 
that $L(z,w)$ has a simple pole in $z$ at $w$ with principal part
$$\frac{1}{2\pi}\frac{1}{(z-w)}.$$
The Szeg\H o and Garabedian kernels have extension,
holomorphicity and antiholomorphicity,
and smoothness properties analogous to those of
$K$ and $\Lambda$, respectively. Finally, $L(z,w)$
has the important property that $L(z,w)\ne0$ if
$z\ne w$ in $\Obar$.

The function $4\pi L(z,w)^2$ is like $\Lambda(z,w)$
in that it has a double pole in $z$ at $w$ with
principal part
$$\frac{1}{\pi}\frac{1}{(z-w)^2}.$$
(The vanishing of the residue term follows from the
fact that
$\int_{b\Om} L(z,w)^2 \,dz$ is equal to minus the
conjugate of
$\int_{b\Om} S(z,w)^2 \,dz$, which is zero by
Cauchy's theorem.)

Standard proofs of identity~\eqref{BS} use the
fact that the one-form
$$(K(z,w)-4\pi S(z,w)^2)\,dz$$
is equal to minus the conjugate of
$$(\Lambda(z,w)-4\pi L(z,w)^2)\,dz$$
on the boundary, which is also a holomorphic one-form
because the poles cancel out, and so the given one-form
extends to the double as a {\it holomorphic\/} one-form
$H_w$. Note that we may write $$H_w=\kappa_w-4\pi\sigma_w,$$
where it is understood that the double poles cancel out.
Such holomorphic one-forms are well-known to be generated
by the $(n-1)$ holomorphic one-forms that are equal to
$F_j'(z)\,dz$ on $\Om$ and equal to minus the conjugate
of $F_j'(z)\,dz$ on the back side, $j=1,\dots,n-1$.
(See \cite[p.~135]{B} for a more elementary proof of
(\ref{BS}).)

Identity~\eqref{BS} shows that the complex linear
span of the functions of $z$ given by
$$K(z,w)-4\pi S(z,w)^2$$
as $w$ ranges over $\Om$ is at most an $n-1$ dimensional
vector space $W$.  We will prove Hejhal's theorem in
the next section by showing that $W$ has to be {\it at
least\/} $n-1$ dimensional because the
$\beta$-periods of linear combinations of $H_w$
as $w$ ranges over $\Om$ can be made to be anything
we like.

The motivation for the
proof in the next section is that the terms $K(z,w)\,dz$ do
not contribute to the value of the $\beta$-periods of
the extension $H_w$ of $(K(z,w)-4\pi S(z,w)^2)\,dz$ to
the double, and the terms $L(z,w)^2\,dz$ can be
used to manipulate the value of the $\beta$-periods to be
anything we like. At the heart of this result is a
density theorem for the Garabedian kernel.
Given a point $a$ in $\Om$, let $L^0(z,a)$ denote the
Garabedian kernel $L(z,a)$ and let $L^m(z,a)$ denote
the derivative $\frac{\dee^m}{\dee w^m}L(z,w)$ evaluated
at $w=a$. Similarly, use a superscript $m$ to indicate
differentiation of the Szeg\H o kernel with respect to
$\bar w$ when $w$ is the second variable in $S(z,w)$.
We will show that the ``Garabedian span at $a$'', which
is the complex linear span of the functions $L^m(z,a)$
as $m$ ranges over the natural numbers, can be used to
approximate functions on the curves $\sigma_j$ that will
lead to elements in the linear span of $H_w$ as
$w$ ranges over $\Om$ with arbitrary $\beta$-periods.
The ``Szeg\H o span at $a$'' is the complex linear span
of the functions $S^m(z,a)$ as $m$ ranges over the natural
numbers.

\section{Proof that $[\lambda_{ij}]$ is nonsingular}
\label{sec3}

We continue to assume that $\Om$ is a bounded domain
bounded by $n>1$ nonintersecting smooth real analytic
curves, and we use the notations and definitions of
the previous section.

The inspiration for the new proof we are about to
give comes from the proof of Lemma~5.1 in
\cite{BGS}, and is yet another reason to view
Mergelyan's theorem as the theorem that is just
too good to be true.

Because the argument needed from Lemma~5.1 of \cite{BGS}
is short after all the machinery we have set up, we
include it here for completeness.
Given a small $\epsilon>0$, let $V$ denote the set
of points in $\C$ that are a distance less than or
equal to $\epsilon$ from $b\Om$. We will shrink
$\epsilon$ as needed in what follows; keep in mind
that $V$ depends on $\epsilon$.
For $j=1,\dots,n-1$, let $\phi_j$ be a continuous
function on the closure of $\sigma_j$ that is equal to
zero on $V\cap\sigma_j$. Thus, $\phi_j$ is zero near
both endpoints of $\sigma_j$. We assume
that $\epsilon$ is small enough that a large open
subset of each curve $\sigma_j$ is not contained in $V$.

One version of Mergelyan's theorem states that, given
a compact set $K$ in the complex plane such that $\C-K$
has only finitely many components and a complex
valued continuous function $\phi$ on $K$ that is
holomorphic in the interior of $K$, there is a rational
function with possible poles only in $\C-K$ that is as
close in the uniform norm as desired to $\phi$ on $K$.
(See Exercise~1 of Chapter~20 in Rudin \cite{R} or Greene
and Krantz \cite[p.~374]{GK}.)

Let $$K=V \cup\left(\cup_{j=1}^{n-1}\sigma_j\right),$$
and let $U=\Om-K$. Note that $U$ is a simply connected
domain contained in $\Om$ if $\epsilon$ is small
enough. By Mergelyan's theorem, there is a rational
function $r(z)$ with possible poles only in $\C-K$
that is as close in the uniform norm as desired to
zero on $V$ and $\phi_j$ on each $\sigma_j$.
As in Stein and Shakarchi \cite[p.~63]{SS}
(and as in many proofs of Runge's theorem) we may slide
the poles of $r(z)$ that fall in $\Om$ to a single point
$a$ in $U\subset\Om$. Let $N$ denote the order of the pole
of $r(z)$ at $a$.

The proof hinges on the following application of
the residue theorem,
\begin{equation}
\label{key}
\frac{2\pi}{2\pi i}\int_{b\Om}r(w)\,L(w,z)\ dw =
r(z)-\sum_{m=0}^{N-1} c_m L^m(z,a)
\end{equation}
for $z\in\Om$ not equal to $a$.
Note that we have used the facts that the principal part
of $L(w,z)$ is
$$\frac{1}{2\pi}\frac{1}{(w-z)}$$
as a function of $w$ and that the only pole of $r(w)$ is
a pole of order $N$ at $a$. The coefficients $c_m$ only
depend on the principal part of $r(z)$ at $a$.
This identity will allow us to approximate $r(z)$ on $K$
by functions in the Garabedian span at $a$. Indeed, using
identity (\ref{SL}) reveals that the left hand side of the
equation is equal to
$$\int_{b\Om}S(z,w)r(w)\ ds_w,$$
where $ds_w$ denotes arc length measure in the $w$-variable,
and this integral is equal to the Szeg\H o projection of
$r(w)$ at the point $z$. Since $r(w)$ can be taken to be
arbitrarily $C^\infty$ close to the zero function on the
boundary, and since the Szeg\H o projection is a continuous
operator from $C^\infty(b\Om)$ to itself (see
\cite[p.~15]{B}), the left member of (\ref{key}) is uniformly small in $z$ 
on $\Obar$. Thus $\mathcal{L}(z)\approx r(z)$ on $\Obar$ for the
approximating element $\mathcal{L}(z)=\sum_{m=0}^{N-1} c_m L^m(z,a)$ in the
Garabedian span at $a$. In particular, $\mathcal{L}(z)\approx 0$
on $b\Omega$, $\mathcal{L}(z)\approx \varphi_j(z)$ on $\sigma_j$.

For a fixed $k$ we now let $\phi_j\equiv0$ for $j\ne k$ and let 
$\phi_k$ be a continuous function on $\sigma_k$ that is zero on
$\sigma_k\cap V$ and such that 
$\int_{\sigma_k} \phi_k(z) L(z,a) \,dz=1$.
(Keep in mind that $L(z,a)$ is nonvanishing on $\Om-\{a\}$.)
We now claim that, for the element $\mathcal{L}(z)$ in the 
Garabedian span constructed above, $\mathcal{L}(z)L(z,a)$
has integrals with respect to $dz$ that are close to zero along
$\sigma_j$, $j\ne k$, and close to one along $\sigma_k$.
Differentiating identity (\ref{SL}) $m$-times with respect to $w$
and then multiplying it by (\ref{SL}) and setting $w$ equal to $a$
reveals that
$$\overline{S^m(z,a)S(z,a)}\,d\bar z=-L^m(z,a)L(z,a)\,dz.$$
Hence, there is
an element $\mathcal{S}(z)$ in the Szeg\H o span at $a$
such that
\begin{equation}
\label{SL3}
\overline{{\mathcal{S}}(z)S(z,a)}\,d\bar z=
{\mathcal{L}}(z)L(z,a)\,dz,
\end{equation}
and because the function
${\mathcal{L}}(z)L(z,a)$
on the right hand side is small
on the boundary of $\Om$, the function
${\mathcal{S}}(z)S(z,a)$ must also be small, and consequently
also small on $\Obar$. Taking the conjugate of (\ref{SL3})
reveals that the one-form
${\mathcal{S}}(z)S(z,a)\,dz$ extends to the double as
a meromorphic one-form $s_a$ by setting it equal to the
conjugate of
${\mathcal{L}}(z)L(z,a)\,dz$ on the back side.
Our construction shows that the $\beta_j$ periods of $s_a$
are close to zero for $j\ne k$ and the $\beta_k$ period is
close to one.

The next step in the proof is to show that there are linear
combinations of the holomorphic one-forms $H_w$ as $w$
ranges over $\Om$ that have $\beta$-periods close to the
$\beta$-periods of the meromorphic one-form $s_a$ constructed
above. Since the $\beta$-periods of $\kappa_w$ are all zero,
it will suffice to find linear combinations of the meromorphic
one-forms $\sigma_w$ that have $\beta$-periods close to the
$\beta$-periods of $s_a$. This turns out to be a rather
elementary exercise due to the following observations.
For small complex $h$,
\begin{align*}
L^1(z,a)L(z,a) &\approx
\frac{L(z,a+h)-L(z,a)}{h}\,\cdot\,L(z,a) \\
&\approx
\frac{L(z,a+h)-L(z,a)}{h}\,\cdot\, \frac{L(z,a+h)+L(z,a)}{2},
\end{align*}
and this last term is a linear combination of the
squares $L(z,a+h)^2$ and $L(z,a)^2$.
(Note that $L^0(z,a)L(z,a)= L(z,a)^2$ is a square; that's why we
skipped it.)
Next,
\begin{align*}
L^2(z,a)L(z,a) &\approx
\frac{L^1(z,a+h)-L^1(z,a)}{h}\ L(z,a) \\
&\approx
\frac{\frac{L(z,a+h+k)-L(z,a+h)}{k}
-\frac{L(z,a+k)-L(z,a)}{k}}{h}\ L(z,a),
\end{align*}
which can be pulled apart into linear combinations of
terms of the form
$$\left[L(z,a+h_1)-L(z,a+h_2)\right]\,\cdot\,L(z,a)$$
that can be approximated by
$$\left[L(z,a+h_1)-L(z,a+h_2)\right]\,\cdot\,
\left(\frac{L(z,a+h_1)+L(z,a+h_2)}{2}\right),$$
which again is a linear combination of squares.
This process can be continued to all higher order terms.

We may now state that there are linear combinations of
the holomorphic one-forms $H_w$ as $w$ ranges over a
small disc $D_\epsilon(a)\subset\Om$
with $\beta$-periods close to any prescribed set of values. 
We must conclude that the linear span is $n-1$
dimensional, and that therefore, the matrix
$[\lambda_{ij}]$ must be nonsingular.
Notice that we above constructed linear combinations
of the one-forms $\sigma_w$ (for $w$ running over $\Omega$)
which have $\beta$-periods essentially equal to those of $s_a$
and hence generating a dense set of period vectors in $\C^{n-1}$.
Therefore Theorem~\ref{thm:main} is now proved.

By showing that $[\lambda_{ij}]$ is nonsingular, we have
proved that the family of functions of $z$ of the form
$$\sum_{i,j=1}^{n-1}
\lambda_{ij}F_i'(z)\overline{F_j'(w)}$$
spans an $n-1$ dimensional vector space of functions
of $z$ as $w$ ranges over any disc
$D_\epsilon(a)\subset\Om$. This implies also that the vectors
$$(F_1'(w),F_2'(w),\dots,F_{n-1}'(w))$$
must span $\C^{n-1}$ as $w$ ranges over
$D_\epsilon(a)$.

\section{Hejhal's theorem in the two-connected case}
\label{sec4}

We now turn to showing directly (without using Theorem~\ref{thm:main}) the positivity of the matrix
$[\lambda_{ij}]$ when $\Om$ is two-connected, i.e.,
that $\lambda_{11}>0$. Since there is only one function
$F_1'$ and one constant $\lambda_{11}$ in \eqref{BL},
we will drop the subscripts.

It was proved in \cite{B1} (see also \cite[p.~149]{B}) that,
for $a$ in one of the boundary curves of $\Om$, the Szeg\H o
kernel $S(a,w)$ has exactly one zero in $\Obar-\{a\}$
in $w$ at a point $b$ in the other boundary curve
of $\Om$. Hence, $\eqref{BS}$ yields that
\begin{equation}
K(a,b)=\lambda F'(a)\,\overline{F'(b)}.
\end{equation}
Multiply this equation by $T(a)\,\overline{T(b)}$ to
obtain
\begin{equation}
\label{KFs}
T(a)K(a,b)\,\overline{T(b)}=\lambda
F'(a)T(a)\,\overline{F'(b)}\,\overline{T(b)}.
\end{equation}

Now the positivity of $\lambda$ follows from the following two consequences of the 
Hopf maximium principle (Hopf lemma):
\begin{equation}
\label{TKT}
T(a)K(a,b)\,\overline{T(b)}<0,
\end{equation}
\begin{equation}
\label{FTFT}
F'(a)T(a)\,\overline{F'(b)}\,\overline{T(b)}< 0.
\end{equation}
In terms of the outward normal derivatives of the Green's function and of the harmonic measure 
$\omega$ these two inequalities express, via (\ref{BLG}), (\ref{Fw}), that
$$
\frac{\partial^2 G(a,b)}{\partial n_a \partial n_b}>0,
$$
$$
\frac{\partial \omega(a)}{\partial n_a}\cdot\frac{\partial \omega(b)}{\partial n_b}<0.
$$
The first inequality actually holds for any two $a,b\in b\Omega$, $a\ne b$, and with
$\Omega$ of arbitrary connectivity. It expresses that the Poisson type kernel 
$p(z,a)=- \frac{1}{2\pi}\frac{\partial G(z,a)}{\partial n_a}$ ($z\in\Om$) attains its minumum value (namely zero) at 
any point on the boundary (for example $z=b$) with a strictly negative slope. 
Similarly, the second inequality
says that $\omega$ has strictly positive (negative) normal derivative on a boundary component on which
it takes its maximum (minimum) value. 

Using that the complex number $-iT(a)$ can be identified with the outward normal vector of $b\Omega$ at $a$
and that $2\dee u/{\dee}\bar{z}$ can be identified with the gradient when $u$ is a real-valued function, 
it follows that if $u$ is constant on $b\Omega$ then 
$$
-2i\frac{\dee u}{\dee z}(a)T(a) \text{ is real and equals } \frac{\dee u}{\dee n_a}.
$$
Using then (\ref{Fw}) and (\ref{BLG}) the inequalities 
(\ref{TKT}), (\ref{FTFT}) follow easily.

We remark that for the Szeg\H o kernel one has
$$
T(a)S(a,b)^2\overline{T(b)}\leq 0,
$$
where equality can be attained (something we have already used). The proof follows on
using $S(a,b)T(a)= i \overline{L(a,b)} $ and $L(a,b)=-L(b,a)$, whereby the inequality becomes
$L(a,b)\overline{L(a,b)}\geq 0$.

\section{Suita's proof that $K(a,a)-4\pi S(a,a)^2>0$} 
\label{sec6}

In the general $n$-connected setting, we have shown
that the matrix $[\lambda_{ij}]$ is nonsingular, and we have 
also proved that it is positive definite in the $2$-connected 
case. We became enthralled with the idea of setting
up an induction via a homotopy argument to deduce
Hejal's complete result that the matrix is positive
definite in general, but we have not been able to
complete the argument. In our quest to find a
shorter, simpler proof of Hejhal's result, we hoped
to use Suita's \cite{Su} beautiful and short proof that
$$K(a,a)-4\pi S(a,a)^2>0$$
as a key step. Since Suita's
result is very much in the spirit of this paper and
since we have set up the tools and notation necessary
to describe it, we include Suita's proof here in case
our readers are inspired to someday complete the plan
of our proof.

Consider the multi-valued function
$$
F(z)=\exp(-G(z,a)-i G^*(z,a))
$$
where $G^*(z,a)$ represents a multi-valued harmonic
conjugate for the Green's function $G(z,a)$.
Note that, because $G(z,a)=-\ln|z-a|+u_a(z)$, where
$u_a(z)$ is the harmonic function that solves the
Dirichlet problem with boundary data $\ln|z-a|$, the
multi-valued function $F$ is bounded in modulus by
a constant times $|z-a|$ near $z=a$. In fact,
$$|F(z)|^2=\exp(-2G(z,a))=|z-a|^2\exp(-2u_a(z))$$
is a single-valued function that is in
$C^\infty(\Obar)$ and is $C^\infty$-smooth
up to the boundary and equal to one there.
The Cauchy-Riemann
equations yield that the complex derivative of an
analytic function $u+iv$ is $u_x-iu_y=2\dee u/\dee z$.
Thus, the complex derivative of the locally defined
analytic function $F(z)$ is given by
$$F'(z)=-2\frac{\dee G(z,a)}{\dee z}
\exp(-G(z,a)-i G^*(z,a))
=-2\frac{\dee G(z,a)}{\dee z}F(z)$$
and so, using the shorthand notation $G=G(z,a)$,
$$|F'(z)|=2\left|\frac{\dee G}{\dee z}\right|
\exp(-G).$$
Since the complex conjugate of $\frac{\dee G}{\dee z}$ is
$\frac{\dee G}{\dee\bar z}$, we may also write
$$|F'(z)|^2=4\frac{\dee G}{\dee z}
\frac{\dee G}{\dee\bar z}
\exp(-2G).$$
Now it is clear that, even though $F$ might be multi-valued,
the quotient $F'/F$ is equal to $-2\frac{\dee G}{\dee z}$ and
is a single-valued analytic function on $\Om-\{a\}$ with a simple
pole at $a$. Using these facts we have collected,
we may compute
$$
\int_\Omega |F'(z)|^2\ dxdy=
4\int_\Omega \frac{\partial G}{\partial z}\frac{\partial G}{\partial \bar{z}}
e^{-2G}\ \left(\frac{1}{2i}d\bar z\wedge d z\right)=
$$ 
$$
=-2i\int_\Omega
\frac{\partial}{\partial\bar{z}}
\left(-\frac{1}{2}e^{-2G}\frac{\partial G}{\partial
z}\right)\,d\bar{z}\wedge dz
=i\int_{b\Omega}e^{-2G}\frac{\partial G}{\partial z}\,dz=
$$
$$
=i\int_{b\Omega}\frac{\partial G}{\partial
z}\,dz=i\int_{b\Om}\left(-\frac{1}{2(z-a)}+{\rm regular}\right)dz=\pi.
$$

Next, let 
$$
f(z)=\frac{S(z,a)}{L(z,a)}
$$
be the Ahlfors map, which is an $n$-to-one branched covering
map of $\Om$ onto the unit disc that is $C^\infty$-smooth up
to the boundary (see \cite[Chap.~13]{B}). It is the solution
of the extremal problem to maximize the modulus
of the derivative at $a$ among all analytic functions
that map $\Om$ into the unit disc, normalized so that
the derivative at $a$ is real and positive. Well known
properties of $f$ include that $f(a)=0$, $f'(a)=2\pi S(a,a)$,
$|f|<1$ on $\Om$, and $|f|=1$ on $b\Om$.

Suita's proof of the inequality is a fiendishly clever
comparison of the Ahlfors map to the mapping $F$,
which can be thought of as a multi-valued substitute
for the Riemann map in the multiply connected setting.
(The Ahlfors map can also be thought of as a
non-one-to-one substitute for the Riemann map in the
multiply connected case.)

Write $f(z)=(z-a)g(z)$ and note that $g(z)$ is an analytic
function on $\Om$ that is nonzero at $a$, has $(n-1)$
zeroes (counted with multiplicity) in $\Om-\{a\}$, and
is $C^\infty$-smooth up to the boundary.
We now have
$$
\ln \left|\frac{f(z)}{F(z)}\right|^2= 2\ln |g(z)|+ 2u_a(z)).
$$
Note that $|f|=|F|=1$ on $b\Omega$. So
$\ln \left|\frac{f(z)}{F(z)}\right|^2$ is a harmonic function
on $\Om$ minus the finitely many zeroes of $g$ in $\Om$,
where it tends to minus infinity, and it
equals zero on the boundary of $\Om$. Since $g$ must
have at least one zero in $\Om$, it follows from the
maximum principle that
$\ln |\frac{f(z)}{F(z)}|<0$ in $\Omega$. Hence
\begin{equation}
\label{fF}
|f(z)|< |F(z)|\quad \text{for } z\in \Omega-\{a\}.
\end{equation}

The combination
$$
h(z):=f(z)\frac{F'(z)}{F(z)}=
\frac{S(z,a)}{L(z,a)}\frac{F'(z)}{F(z)}=
-2 \frac{S(z,a)}{L(z,a)}\,\frac{\partial G(z,a)}{\partial z}
$$
is analytic in $\Omega$ if we set $h(a)=f'(a)=2\pi S(a,a)$.
Therefore
$$
2\pi S(a,a)=h(a)=\int_\Omega h(z)\overline{K(z,a)}\,dxdy,
$$
and so, using (\ref{fF}),
$$
4\pi^2S(a,a)^2\leq \int_\Omega |h(z) ^2 |\,dxdy\cdot
\int_\Omega | K(z,a)|^2\,dxdy
$$
$$
=\int_\Omega \left|\frac{f(z)}{F(z)}\right|^2 | F'(z)|^2\,dxdy\cdot K(a,a)
$$
$$
<\int_\Omega | F'(z)|^2dxdy\cdot K(a,a) =\pi K(a,a).
$$
Now the desired inequality is proved. 

\section{The three-connected case}
\label{sec5}




Suita's result 
allows us to deduce Hejhal's theorem from rather basic facts in the case $n=3$.
The matrix $[\lambda_{ij}]$ is easily seen to be
hermitian (it is in fact real and symmetric). Hence it is diagonalizable and we may
introduce a new basis $U_k'(z)$ for the linear span
of the functions $F_j'$ which are linear combinations
of the $F_j'$ with real coefficients such that
$$K(z,w)-4\pi S(z,w)^2=\sum_{i=1}^{2}\mu_i U_i'(z)\overline{U_i'(w)}.$$
We know that the $\mu_i$ (the eigenvalues of $[\lambda_{ij}]$)
are real and nonzero (because $[\lambda_{ij}]$ is nonsingular).
We now consider the zeroes of the $U_k'$.
Since $U_k'(z)T(z)=-\overline{U_k'(z)T(z)}$ on the boundary, 
the generalized argument principle, that allows zeroes
on the boundary that are counted with a factor of
one-half in front, shows that each $U_k'(z)$ has either one zero in $\Omega$
or two zeroes on the boundary.

Now, Suita's result yields that
$$\mu_1|U_1'(a)|^2+\mu_2|U_2'(a)|^2 > 0$$
for any point $a$ in $\Om$. 
Choosing $a$ to be the zero of $U_2'$, in case that
zero is in $\Om$, shows that $\mu_1$ is necessarily 
positive. If instead $U_2'$ has two zeroes on the boundary
it still follows, by sliding such a zero a little into $\Omega$,
that $\mu_1$ cannot  be strictly negative. This is because
$U_1'(a)\ne 0$ at any $a$ with $U_2'(a)=0$ (not all holomorphic
differentials in a basis can vanish at the same point). Invoking
Theorem~\ref{thm:main} we then actually have $\mu_1>0$
again. Similarly $\mu_2>0$. This completes the proof.

\section{Wishful thinking}
\label{sec7}

We have feelings and urges about Hejhal's theorem
that we'd like to share here. Hejhal's proof of the
positive definiteness of the lambda matrix was long and
technical and used theta functions in a key way.
We would like to come up with alternative ways
of understanding the result, ways that might be
quicker and more elementary. One potentially easier
way to understand the proof might involve a homotopy
argument. As a smoothly bounded $n$-connected
domain varies in a $C^\infty$ way, the kernel functions
and lambda coefficients vary in a $C^\infty$ way, too.
Because we have shown that the matrix is nonsingular, 
an eigenvalue of the matrix could not pass through
zero to change sign under smooth perturbations. We
have shown the $1\times1$ matrix in the $2$-connected
case is positive definite. If we could show that
there is just {\it one\/} $n$-connected domain for
which the matrix is positive definite, then all
$n$-connected domains must share that property
since they are all smoothly homotopic. Our idea is
to let one of the holes in an $n$-connected domain
shrink down to a point, perhaps heading off to the
outer boundary curve as a shrinking circle. If one
understood the asymptotic behavior of the kernel
functions in the limit under this process, then
perhaps one could deduce the positive definiteness
of an $n$-connected domain near the limiting domain
from knowing the result in the $(n-1)$-connected case.
The $2$-connected case would start the induction
process. In fact, understanding the asymptotic
behavior of the Szeg\H o kernel is all that would
be needed, as we now explain.

Let $dF_j$ denote the holomorphic $1$-form that
is $F_j'(z)dz$ on $\Om$ and $-\overline{F_j'(z)dz}$
on the backside of the double. Because
$F_j'=2\dee\omega_j/dz$, it is easy to check that
the $\beta$-periods of $dF_j$ are given by two times
the delta function, i.e.,
$$\int_{\beta_k}dF_j=2\delta_{kj}.$$
Let $\kappa_w$ and $\sigma_w$ denote the $1$-forms
that we introduced in \S\ref{sec2}. Recall that
the $\beta$-periods of $\kappa_w$ are zero. Hence,
when we take the $\beta$-periods of the identity
(\ref{BS}), we get
a formula that only involves the Szeg\H o kernel,
the $\lambda$ coefficients, and the function $F_j'$.
Taking this idea further, let $\kappa$ denote the
$(1,1)$-form gotten from extending $dz\,K(z,w)\,d\bar w$
to the double cross the double and let
$\sigma$ denote the $(1,1)$-form gotten by
extending $dz\,S(z,w)^2\,d\bar w$. Our observations
about the $\beta$-periods reveal the known 
(compare equation (30), p.231, in \cite{HS}) formula
$$4\pi\int_{\beta_i}\int_{\beta_j}\sigma dz\,d\bar w=
-4\lambda_{ij}$$
when $i\ne j$. For the case $i=j$, we let $\widetilde{\beta_i}$
denote a curve obtained from sliding $\beta_i$ along the
$i$-th handle some distance. In this case, we have
$$4\pi\int_{\beta_i}\int_{\widetilde{\beta_i}}\sigma dz\,d\bar w=
-4\lambda_{ii}.$$
This shows, strangely enough, that the $\lambda$ matrix
only depends on the Szeg\H o kernel! Hence, if we
understood the asymptotic behavior of the Szeg\H o
kernel as a circular hole shrinks away to nothing
as it heads toward the outer boundary, we might be
able to reduce Hejhal's complete result in the $n$-connected
case to the $(n-1)$-connected case to complete our
wished for induction.

Using the Kerzman-Stein integral formula for the
Szeg\H o kernel as in \cite[p.~153-157]{B}
seemed like a promising way to deduce that the
kernel functions of the shrinking hole domains
converge in a strong sense to the kernel functions
of the lower-connectivity limit domain. Assuming
the shrinking hole is circular simplifies some
of the arguments because the Kerzman-Stein kernel
is zero for $(z,w)$ on the same circular boundary curve. 
A different idea would be to use slit models for domains $\Omega$
having a hyperelliptic double and use explicit formulas, due to
Barker \cite{Ba}, for the Szeg\H o  kernel for such domains.
But alas, the needed results eluded us.
We happily continue to ruminate upon it.

\medskip

\noindent
{\bf Data availability statement:} This research does
not depend on data.

\end{document}